\documentstyle[12pt]{article}

  \newcommand{\const}{\rm const}
  
  \newcommand{\supp}{\rm supp}
  
  \newcommand{\sign}{\rm sign}
  
  \newcommand{\vraisup}{\rm vraisup}

  \newcommand{\argmin}{\rm argmin}
  
  \newcommand{\mod}{\rm  mod}

\begin{document}

 \begin{center}

 \ {\bf Grand Lebesgue Spaces norm estimates } \par

\vspace{4mm}

{\bf for multivariate functional operations}\par

\vspace{5mm}

 \ {\bf   Ostrovsky E., Sirota L.}\\

\end{center}

\vspace{4mm}

 ISRAEL, Ramat Gan, Bar-Ilan University, \\

 department  of Mathematic and Statistics, 59200. \\

\vspace{4mm}

e-mails: \  eugostrovsky@list.ru \\
sirota3@bezeqint.net \\

\begin{center}

\vspace{4mm}

  \ {\bf Abstract} \\

\end{center}

\vspace{4mm}

  \ We intend to derive the moment and exponential tail estimates for the so-called bivariate or more generally multivariate {\it functional} operations, not necessary
 to be linear or even multilinear. We will show also the strong or at last weak (i.e. up to multiplicative constant) exactness of obtained estimates. \par

 \vspace{4mm}

\begin{center}

\vspace{4mm}

  \ {\bf Key words and phrases.} \\

\end{center}

 \ Binary and multivariate operators and operations; Banach, Lebesgue-Riesz, Orlicz and Grand Lebesgue Spaces (GLS); moment and tail estimates, ordinary convolution and
 infimal convolution, examples, exact or weak exact value of constants, upper and lower ordinary and exponential estimates, generating function, tail function, singular integral operators,
 multivariate multiplicative, tensor, Haussdorf, Hilbert, maximal, pseudo-differential, Hardy-Littlewood and other operations. \par

 \vspace{6mm}

 \section{ Introduction. Statement of problem.}

 \vspace{4mm}

 \ Let $ \ (X,F,\mu), \ (X_i, F_i, \mu_i), \ i = 1,2,\ldots,d  \ $ be measurable spaces equipped with
 non - trivial measures  $ \ \mu, \mu_i, \ $  not necessary to be probabilistic or bounded.
 We denote by $ \  L(p, x_i) \ $ the classical Lebesgue-Riesz spaces of measurable functions $ \ f_i: X_i \to R \ $  having a finite norm

 $$
 |f_i|L(p,X_i) = |f_i|_p := \left[ \ \int_{X_i} |f_i(x_i)|^{p} \ \mu_i(dx_i) \   \right]^{1/p}, \ p \in [1,\infty). \eqno(1.0)
 $$

\vspace{4mm}

 \ Let also $ \  V: \otimes_{i=1}^d L(q_i,  X_i) \to L(p,X) \ $ be multivariate functional {\it operation} which maps the tensor product
$   \ \otimes_{i=1}^d L(q_i,  X_i), \  p_i \in [1,\infty) \  $  into the Lebesgue-Riesz space $ \ L(p,X): \ $

$$
g = g(x) =  V[f_1, f_2,   \ldots, f_d](x) \in L(p,X). \eqno(1.1)
$$

 \ It will be presumed more precisely  that there exists  a non-trivial set $ \  D \subset R^d_+  \  $ such that

$$
\forall \ \vec{q} \in D \ \Rightarrow g(\cdot) \in L(p,X).
$$

 \ Let $ \ \tau_i = \tau_i(q_i)\ $ be some continuous inside its domain of definition functions and strictly monotonically increasing functions, for instance

$$
\tau_i(q_i) =  \beta_i (q_i)^{\gamma_i}, \ \beta_i, \gamma_i = \const \in ( 0, \infty), \ q_i \in [1,\infty).
$$

 \ Denote $ \ \vec{q} := (q_1,q_2, \ldots, q_d) = \{  q_i \},  \ q_i \in [1,\infty) \ $ and define analogously a vector
 $ \ \vec{\alpha} := \{\alpha(i) \} = (\alpha_1, \alpha_2, \ldots, \alpha_d), \ \alpha_i \ge 0;  \ $

$$
\vec{q}^{ \vec{\alpha}  } \stackrel{def}{=} \prod_{i=1}^d \left(q_i \right)^{\alpha_i}, \
| \ \vec{f} \ |_{\vec{q}}^{ \vec{\alpha} } \stackrel{def}{=} \prod_{i=1}^d \left(|f_i|_{q_i} \right)^{\alpha_i},
$$
or more generally

$$
\vec{\tau}(\vec{q}) \stackrel{def}{=} \left\{ \tau_i(q_i) \ \right\}, \  \vec{\tau}(\vec{q})^{ \vec{\alpha} } \ \stackrel{def}{=} \prod_{i=1}^d \left[ \ \tau_i(q_i) \ \right]^{\alpha_i},
$$

$$
| \ \vec{f} \ |_{\vec{\tau}(\vec{q})}^{ \vec{\alpha} }  \stackrel{def}{=} \prod_{i=1}^d \left(|f_i|_{\tau_i(q_i)} \right)^{\alpha_i},
$$

 \ We assume in the sequel that this operation satisfies the following condition. \par

\vspace{4mm}

 \ {\bf Condition 1.1.} \ \ {\it  There exists a non-trivial domain $ \ D \subset \otimes_{i=1}^d [1, \ \infty) \subset  R^d_+ \ $ and
 a certain function $ \Theta = \Theta(\vec{q}) \in [1,\infty), \  \vec{q} \in D  \ $ such that for the value} $ \ p = \Theta(\vec{q}) \in [1, \ \infty) \ $
{\it the following  inequality holds true}

$$
|g|_p = |V(f_1, f_2, \ldots, f_d)|_{\Theta(\vec{q})} \le K(\vec{q}) \cdot  | \ \vec{f} \ |_{\vec{\tau}(\vec{q})}^{ \vec{\alpha} }, \
$$
where
$$
K(\vec{q}) = K_{\alpha,\beta, \gamma, \tau}(\vec{q}) = \const < \infty \ \ \Longleftrightarrow \vec{q} \in D. \eqno(1.2)
$$

  \  One can extend formally the definition this function as follows:

$$
 K_{\alpha,\beta, \gamma, \tau}(\vec{q})  := + \infty, \ \vec{q} \notin D.
$$

\vspace{4mm}

\ As for the function $ \ K(\vec{q}): \ $  it may consists on the factors of the form for instance

$$
K_{1,i}(q_i) :=  \beta_i (q_i)^{\gamma_i}, \ \beta_i, \gamma_i = \const \in ( 0, \infty), \ q_i \in [1,\infty);
$$
or

$$
K_{2,j}(q_j) := \frac{ \beta_j (q_j)^{\gamma_j}}{ (q_j - 1)^{\delta_j}   }, \
\beta_j, \gamma_j, \delta_j = \const \in ( 0, \infty), \ q_j \in (1,\infty);
$$
or at last

$$
K_{3,l}(q_l) = C \ (q_l - a_l)^{-c_l} \ (b_l -  q_l)^{-s_l}, \ q_l \in (a_l, \ b_l),  \
$$

$$
 1 \le a_l < b_l < \infty; \ c_l, \ s_l = \const \ge 0;
$$
and so one. \par

 \ An example:

$$
K(q_1,q_2,q_3) = C (q_1)^{\gamma_1} \cdot \frac{(q_2)^{\gamma_2}}{(q_2 - 1)^{\delta_2}} \cdot (q_3 - a_3)^{-c_3}  \ (b_3 - q_3)^{-s_3},
$$
where

$$
C \in (0,\infty), \ \gamma_1,\delta_2, c_3, s_3 > 0,  \ \gamma_2 \ge 0, \ 1 \le a_3 < b_3 < \infty,
$$
and

$$
q_1 \in [1,\infty); \  q_2 \in (1,\infty); \ q_2 \in (a_3, b_3).
$$

\vspace{4mm}

 \ {\it One can  choose as the  function } $ \ K(\vec{q}) \ $  {\it its minimal value:}

$$
|g|_p = |V(f_1, f_2, \ldots, f_d)|_{\Theta(\vec{q})} \le \overline{K}(\vec{q}) \cdot  | \ \vec{f} \ |_{\vec{\tau}(\vec{q})}^{ \vec{\alpha} }, \eqno(1.3)
$$
so that
$$
\overline{K}(\vec{q}) = \overline{K}_{\alpha,\beta, \gamma, \tau}(\vec{q}):=\sup_{\vec{f} \ne 0} \
\sup_{\vec{q} \in D} \left\{ \ \frac{|V(f_1, f_2, \ldots, f_d)|_{\Theta(\vec{q})} }{ \ | \ \vec{f} \ |_{\vec{\tau}(\vec{q})}^{ \vec{\alpha} }} \ \right\}. \eqno(1.3a)
$$

\vspace{4mm}

 \ {\bf  Our goal in this preprint  is a generalization of the estimate (1.2) into the more general spaces, namely, into a so-called Grand Lebesgue Spaces (GLS).  }\par

\vspace{5mm}

{\sc Hereafter we will denote by $ \ c_k = c_k(·), \ C_k = C_k(·), k = 1, 2, . . . , \ $ with or without subscript, some positive finite non-essentially
"constructive" constants, non necessarily at the same at each appearance.} \par

 \ {\sc We will denote also by the symbols $ \  K_j = K_j(d, n, p, p_i q, . . .) \ $ essentially
positive finite functions depending only on the variables} $ \ d, n, p, q, . . \ $ , \par

\vspace{5mm}

 \ Let us bring  some examples. \par

 \vspace{4mm}

  \ {\bf Example 1. Multiplicative bilinear operator.  } \par

\vspace{4mm}

 \ In this example $ \  X = X_1 = X_2, \ \mu = \mu_1 = \mu_2   \ $ and

$$
g(x) = f_1(x) \cdot f_2(x). \eqno(1.4)
$$

 \ One can apply the classical H\"older's inequality

$$
|g|_p \le |f_1|_{\alpha  \ p}  \cdot |f_2|_{\beta \ p}, \ p \ge 1, \eqno(1.4a)
$$
where
$$
\alpha,\beta > 1, \ 1/\alpha + 1/\beta = 1.
$$
 \ Therefore

$$
|g|_p \le \inf_{\alpha, \ \beta} \left[ \ |f_1|_{\alpha \ p} \ |f_2|_{ \beta \ p} \ \right], \ p \ge 1, \eqno(1.4b)
$$
where $ \  \inf \ $ is calculated over all the values $ \ \alpha, \ \beta > 1, \ $ and such that \\
 $ \ 1/\alpha + 1/\beta = 1. \ $ \par

\vspace{4mm}

\ {\bf Example 2. Tensor product.} \par

\vspace{4mm}

 \ Here $ \ (X,F,\mu) = (X_1,F_1,\mu_1) \otimes (X_2,F_2,\mu_2)  \  $ and

$$
g = g(x_1,x_2) = f_1(x_1) \cdot f_2(x_2). \eqno(1.5)
$$
 \ On the other words, both  the  cofactors  $ \ f_1, \ f_2 \ $ are independent in the probabilistic sense. \par

 \ We conclude

$$
|g|_p = |f_1|_p \cdot |f_2|_p, \eqno(1.5a)
$$
but the last relation is true even for all the non-negative values $ \ p; \ p \ge 0. \ $ \par

\vspace{4mm}

{\bf Example 3.  Integral bilinear operator.} \par

\vspace{4mm}

 \ Let us consider now the following integral bilinear (regular) operator

$$
g(x) := V_L[f_1,f_2](x) \stackrel{def}{=} \int_{X_1} \int_{X_2} L(x,x_1,x_2) \ f_1(x_1) \ f_2(x_2) \ \mu_1(dx_1) \ \mu_2(dx_2), \eqno(1.6)
$$
 $ \ p, p_1,\ p_2 \in (1,\infty). \ $
 \ Put as ordinary  $ \ p' = p/(p-1), \  p_j' = p_j/(p_j - 1), \ j = 1,2. \ $
 Denote also  by $ \  l(p, p_1, p_2) = l[L](p, p_1, p_2)  \ $ the following {\it mixed}, or equally {\it anisotropic} norm of the kernel $ \ L(\cdot): \ $

$$
l[L](p, p_1, p_2) := | \ | \ | \ L \ |_{p_2, X_2} \ |_{p_1, X_1}   \ |_{p, X} \stackrel{def}{=}
$$

$$
\left\{ \ \int_X \ \left[ \ \int_{X_2} \left( \ \int_{X_1} \left|L(x, x_1, x_2) \right|^{p_1'} \mu_1(dx_1) \  \right)^{p_2'/p'_1} \ \mu_2(dx_2) \  \right]^{ p/p'_2} \ \mu(dx) \ \right\}^{1/p}.
$$

 \  This notion was introduced  at first by Benedek A. and Panzone R. [3]; see a detail investigation and applications in [5], [31], [32] etc.\par
 \ It is not hard to obtain by means of H\"older's inequality

$$
|g|_p \le l[L](p, p_1, p_2) \cdot |f_1|_{p_1} \ |f_2|_{p_2}.
$$
 \ On the other words, in this example

$$
K_L(p, p_1, p_2) = l[L](p, p_1, p_2),
$$
of course, for all the values of the parameters $ \ (p, \ p_1, \ p_2) \ $  for which the right - hand side is finite. \par

\vspace{4mm}

{\bf Example 4. Classical convolution.} \par

\vspace{4mm}

 \ Let $ \ X, X_i = R^n, i = 1,2; d = 2 \ $ and $ \mu, \ \mu_i \ $ be as before Lebesgue measures. Consider a classical convolution operation

$$
g = f_1*f_2, \ \Longleftrightarrow  g(x) = \int_{R^n} f_1(x-y) \ f_2(y) \ dy. \eqno(1.7)
$$
 \ Let also  the values $ \   p_1, \ p_2 = \const > 1,  \  $ and $ r = \const \in (1,\infty) \ $ are such that

$$
 1 + \frac{1}{r}  =   \frac{1}{p_1} +  \frac{1}{p_2}. \eqno(1.7a)
$$

 \ It will be presumed of course that (here and in the sequel) in this operation

$$
\forall (p_1, \ p_2) \in D \ \Rightarrow  1 \le \frac{1}{p_1} +  \frac{1}{p_2} \le 2.
$$

\vspace{4mm}

 \ The classical Young's inequality tell us that

$$
|g|_r \le G(r,p_1,p_2) \ |f_1|_{p_1} \ |f_2|_{p_2}, \  G(r,p_1,p_2) = \const \le 1.
$$
 \ The exact value of the "constant"  $ \ G(r,p_1,p_2) \ $ was obtained at first by W.Beckner [2]; see also
 H.J.Brascamp and E.H.Lieb [6]:

$$
 G(r,p_1,p_2) = \left[ \ \frac{v(p_1) \ v(p_2)}{v(r)}  \ \right]^n, \eqno(1.7b)
$$
where

$$
v(p) := \left[ \ p^{1/p} \ (p')^{-1/p'} \ \right]^{1/2}, \ \ p' = p/(p-1).
$$

 \ This fact may be easily generalized onto the unimodular local compact topological group. Let us consider the case $ \  X = X_1 = X_2 = [0, 2 \pi]^d \ $
with ordinary Lebesgue measure. The convolution may be defined alike:

$$
g(x) = f_1 * f_2(x) := (2 \ \pi)^{-d} \int_X f_1(x-y) \ f_2(y) \ dy, \eqno(1.7c)
$$
where all the algebraic operations in (1.7b) are understood $ \ \mod(2 \ \pi). \ $  The estimates (1.7a) and (1.7b)  remains true. \par

\vspace{4mm}

{\bf Example 5. Infimal convolution.} \par

\vspace{4mm}

 \ Here $ \  x,y \in R^d;  \ \mu,\mu_1, \mu_2 \  $ are usually Lebesgue measures.

 $$
 g = f_1 \Box f_2 \ \Longleftrightarrow g(x) \stackrel{def}{=} \inf_{y \in R^d} (f_1(x-y) + f_2(y)). \eqno(1.8)
 $$
 \ This operation appears in the theory of optimization, convex analysis etc. \par
 \ Let $ \ p \ $ be arbitrary number  from the set $ \ [1,\infty); \ $ define the value

$$
K(d,p) \stackrel{def}{=} \sup_{ |f_1|_p + |f_2|_p \in (0,\infty)} \left\{ \ \frac{ \ | f_1 \Box f_2|_p \ }{|f_1|_p + |f_2|_p} \  \right\}. \eqno(1.8a)
$$

\vspace{4mm}

{\bf Theorem 1.}

$$
K(d,p) = 2^{d/p}. \eqno(1.8b)
$$

\vspace{4mm}

{\bf Proof. I. Upper bound.} \par

 \vspace{4mm}

 \ We can and will suppose without loss of generality that both the function $ \ f_1 \ $ and $ \ f_2 \ $ are non-negative, as well as the "common" one
 $ \ g. \ $ We derive choosing $ \ y := x/2: $

$$
g(x) \le f_1(x/2) +  f_2(x/2).
$$

 \ Define as ordinary the so - called  dilation operator

$$
T_{\lambda}[f](x) \stackrel{def}{=} f(\lambda \ x), \ \lambda = \const > 0;
$$
then

$$
| T_{\lambda}[f](\cdot)|_p = \left[ \ \int_{R^d} |f(\lambda \ x)|^p \ dx  \ \right]^{1/p} = \left[ \ \int_{R^d} |f(z)|^p \ \lambda^{-d} \ dz    \ \right]^{1/p} = \lambda^{-d/p} \ |f|_p.
$$
 \ We have applying the last relation for the value $ \ \lambda = 1/2: \ $

$$
|g|_p \le 2^{d/p} \left[ \ |f_1|_p  + |f_2|_p \ \right],
$$
therefore $ \ K(d,p) \le 2^{d/p}. \ $ \par

\vspace{4mm}

{\bf II. Lower estimate.} We choose $ \ f_1 = f_2 =:f, \ $ where $ f = f(x), \ x \in R, \ $ is certain  non-negative  smooth with smooth
strictly increasing on some finite non-trivial interval derivative  function. We conclude

$$
t := \argmin_y [\left[ f(x-y) + f(y)   \right] = x/2;
$$
then

$$
g(x) = [f \Box f](x)  = 2 f(x/2),
$$
following

$$
|g|_p =  2 \left|T_{0.5} \ f(\cdot) \right|_p =2 \cdot 2^{d/p} \ |f|_p = 2^{d/p} \left[ |f_1|_p + |f_2|_p   \right].
$$
 \ Thus, $ \ K(d,p) \ge 2^{d/p}. \ $  \par

\ The {\it multidimensional case}   $ \ d \ge 2 \ $ may explored analogously: $ \ f_0(x_1, x_2) := f(x_1) \cdot f(x_2), \ d = 2. \ $  \par

\vspace{4mm}

\ The {\it multivariate case}

$$
 g_m(x) :=  [ \ f_1 \Box f_2 \Box \ldots f_m  \ ](x)\ = [ \ \Box_{j=1}^m f_j \ ](x)
$$
may be investigated quite analogously:

$$
|g_m|_p  \le m^{d/p} \ \sum_{j=1}^m |f_i|_p, \ p \ge 1,
$$
herewith the constant $ \ m^{d/p}  \ $ is the best possible. \par
 \ As far as we know, see for example a recent review of  Thomas Str\"omberg \  [38], this proposition is new. \par

\vspace{4mm}

{\bf Example 6. Pseudo-differential product.} \par

\vspace{4mm}

 \  A so - called pseudo-differential product (PDP) may be defined as follows. $ \ X = X_1 = X_2 = R, \ $ and as above   $ \mu(dx) = \mu_1(dx) = \mu_2(dx) = dx; \ $

$$
g(x) = PD[f_1, \ f_2](x) := \int\ \int_R e^{ i x (\alpha + \beta)  }  \ \sigma(x,\alpha,\beta) \ \tilde{f}_1(\alpha) \ \tilde{f}_2(\beta) \ d \alpha \ d \beta, \eqno(1.9)
$$
where the notation $ \ \tilde{f}(\cdot) \ $ stands for the Fourier transform:

$$
\tilde{f}(\cdot)(\gamma) \stackrel{def}{=} \int_R e^{i \ \gamma \ y}   \ f(y) \ dy.
$$
 \ We impose  on the "kernel -  symbol" function $ \ \sigma = \sigma(x,\alpha,\beta) \ $ the classical
 H\"ormander's condition [4].\par

 \ Here

$$
K(p,p_1,p_2) < \infty \ \Longleftrightarrow \frac{1}{p} = \frac{1}{p_1} + \frac{1}{p_2} < \frac{3}{2}, \ 1 < p,p_1,p_2 \le \infty,  \eqno(1.9a)
$$
see [4] etc.\par
 \ The case of another pseudo  - differential  bilinear operators, for instance, Bochner-Riesz average, may be found in [28], [29].

\vspace{4mm}

{\bf Example 7. Hilbert's bilinear operator.} \par

\vspace{4mm}

 \ In this example under at the same restrictions as in  previous considerations

$$
g(x) = H_{\lambda(1), \ \lambda(2)}[f_1, f_2](x) := v.p. \int_R f_1(x - \lambda_1 \ y)  \ f_2(x - \lambda_2 \ y) \ dy/y,  \eqno(1.10)
$$
$ \ \lambda_1, \ \lambda_2 = \const \in R. \ $ \par
 \ This Hilbert's bilinear operator is the particular case of the last one, namely

$$
\sigma(x,\alpha,\beta) = i \ \pi \ \sign(\lambda_1 \alpha + \lambda_2 \ \beta).
$$

\vspace{4mm}

{\bf Example 8. Maximal operator.} \par

\vspace{4mm}

 \ Define the following {\it maximal multilinear} operator

$$
g(x) = M_R[f_1,f_2,\ldots,f_d](x) \stackrel{def}{=} \sup_{R: x \in R} \prod_{i=1}^d \left\{ \ \frac{1}{|R|} \int_R |f_i(y_i)| \ dy_i \   \right\},  \eqno(1.11)
$$
where $ \ x \in   R^d \ $  and $ \ R \ $ denotes the family of all rectangles in $ \  R^d \ $ with sides parallel to the axes. \par
 \ Here $ \ p_i \in (1,\infty), \ $

$$
\frac{1}{p} = \sum_{i=1}^d  \frac{1}{p_i}, \eqno(1.11a)
$$

$$
K_R(p,p_1,p_2,\ldots,p_d) := C^d \prod_{i=1}^d \frac{p_i}{p_i - 1}, \eqno(1.11b)
$$
so that

$$
|g(\cdot)|_p = |M_R[f_1,f_2,\ldots,f_d](\cdot)|_p \le K_R(p,p_1,p_2,\ldots,p_d) \ \prod_{i=1}^d |f_i|_{p_i}, \eqno(1.11c)
$$
see [28], [29]; see also the reference therein. \par

\vspace{4mm}

{\bf Example 9. Hausdorff's  operation.} \par

\vspace{4mm}

 \ The operator of a form

 $$
g(x) = H_{\Phi, \vec{A} } [\vec{f}](x) := \int_{R^n} \frac{\Phi(t)}{|t|^n} \ \prod_{i=1}^m f_i(A_i(t)) \ dt, \ x \in R^n\eqno(1.12)
 $$
is said to be  {\it Hausdorff  operator.} \par

 \ It is bounded under some natural conditions, see  [12], [13], [7];
 and as before

$$
K_H(p, \vec{p}) < \infty \ \Longleftrightarrow \ \frac{1}{p} = \sum_{i=1}^m  \frac{1}{p_i}, \ p,p_1 \in (1, \ \infty).
$$

 \ More precisely, under these conditions

$$
K_H(p, \vec{p}) \le C(\Phi, \vec{A}, m,n) \cdot \prod_{j=1}^m \frac{p_j^2}{p_j - 1}. \eqno(1.12a)
$$

\vspace{4mm}

{\bf Example 10.} {\it Bounded multiplicative Toeplitz operators on sequence spaces.}

\vspace{4mm}

 \  The  so - called Toeplitz operator acting on the numerical infinite sequences $ \ \{  x_1,x_1, \ldots   \} $ and in the values also in ones by the formulae

$$
g_n  = M[f](x) \stackrel{def}{=} \sum_{k=1}^{\infty} f\left( \frac{n}{k} \right) \ x_k. \eqno(1.13)
$$
 \ Here the function $ \ f(\cdot) \ $ is defined on the set of all positive rational numbers $ \ Q_+, \ $ equipped with usually uniform measure: $ \ \mu_1 \{l/m\} = 1; \ $
the correspondent measures $ \ \mu; \ \mu_2 \ $ defined  on the set of all positive natural numbers is ordinary countable measure with unit value of each number. \par

 \ Nicola Thorn in the  recent article  [39]  has proved the following bilateral estimate

$$
|y|_p \le "1" \cdot |f|_{p_2} \cdot |x|_{p_1}, \eqno(1.13a)
$$
iff

$$
\frac{1}{p }  = 1 - \frac{1}{p_1 } -   \frac{1}{p_2 }, \ p, p_1, \ p_2 \in (1,\infty). \eqno(1.13b)
$$
 \ Herewith the constant "1" in the Thorn estimate is in general case the best possible \ [39]. \par

\vspace{4mm}

 \ Alike  assertions holds true for the multilinear fractional operator, Bochner-Riesz averages, multiple Riesz transform and so one. \par

 \vspace{4mm}

 \section{ Grand Lebesgue Spaces.}

 \vspace{4mm}

 \ We recall here for reader convenience some used further facts about the so-called Grand Lebesgue Spaces (GLS); more information about this GLS
 may be found in articles and monographs [35],  [25], [26], [8], [9],  [10], [15], [16], [20], [22], [23], [34], and so one. \par

\vspace{4mm}

 \ Let   $ \   ( Z, B,  \nu)  \ $ be certain measure space with some non - trivial measure $ \ \nu. $
 \ Let also  $  \psi = \psi(p), \ p \in(a, b), \  \exists a \ge 1, \ \exists b = \const \in (a, \infty]  $ (or   $ \ p \in [a,b) \ $ ) be  bounded
from below:  $  \ \inf \psi(p)  > 0 $ continuous inside the  {\it  semi - open} interval $   \ p \in (a, b) \ $ numerical valued function
such that the auxiliary  function

$$
h(p) = h[\psi](p) \stackrel{def}{=} p \ \ln \psi(p) \eqno(2.0)
$$
is convex. The set of all such a functions will be denoted by $ \ \Psi; \ \Psi = \cup_{a,b: \ 1 \le a < b < \infty}\Psi(a,b). \ $\par

\vspace{4mm}

 \ As ordinary, in this section  for arbitrary measurable function $ \  f: Z \to R \  $

$$
|f|_p = |f|_{p,\nu}  \stackrel{def}{=} \left[ \  \int_Z |f(z)|^p \ \nu(dz) \   \right]^{1/p}, \ p \in [1, \ \infty).
$$

 \ An important example. Let  $ \ \eta \ $ be a measurable function such that there exists $ \  b = \const > 1   \ $  so that
$ \  |\xi|_b < \infty.  \ $  The so-called {\it  natural } $ \ G\Psi_{\eta} \ $  function $ \  \psi_{\eta} = \psi^{(\eta)}(p)  \  $  for the r.v. $ \ \eta \ $
is defined by a formula

$$
\psi^{(\eta)}(p) \stackrel{def}{=} |\eta|_p.
$$

 \ Then $ \  \eta \in G\psi_{\eta} \ $ and

$$
||\eta||G\psi_{\eta} = 1.
$$

\vspace{4mm}

 \ We can and will suppose $ \ a = \inf \{p, \ \psi(p) < \infty \ $ and  correspondingly
$   \ b = \sup \{p, \psi(p) < \infty\},  \ $ so that  $ \  \supp \ \psi = [a, b) \  $  or $ \ \supp \ \psi = [a, b] \ $   or $ \ \supp \ \psi = (a, b] \ $
or at last  $ \ \supp \ \psi = (a, b). \ $  The set of all such a
functions will be denoted by  $ \ \Psi(a,b) = \{  \psi(\cdot)  \}; \ \Psi := \Psi(1,\infty).  $\par

\vspace{4mm}

 \  {\it By definition, the (Banach) Grand Lebesgue Space  \ (GLS)  \ space   $  \ G\psi = G\psi(a,b)  $ consists on all the
numerical  valued  (real or complex)
measurable functions $ \zeta  $ defined on our measurable space $ \ Z = (Z,\ B, \nu) \ $  and having a finite norm}

$$
||\zeta|| = ||\zeta||G\psi \stackrel{def}{=} \sup_{p \in (a,b)} \left\{ \frac{|\zeta|_p}{\psi(p)} \right\}. \eqno(2.1)
$$

\vspace{4mm}

\  The function $ \  \psi =\psi(p) \  $ is named as a {\it  generating function } for this Grand Lebesgue Spaces. \par

\vspace{4mm}

 \ These spaces  are Banach functional space, are complete, and rearrangement invariant in the classical sense,
 and were  investigated in particular in  many  works, see the aforementioned works. \par

\  We refer here  some  used in the sequel facts about these spaces and supplement more. \par

\vspace{4mm}

 \ Define as usually for any measurable function $ \ \zeta: Z \to R \ $ its tail function

 $$
 T_{\zeta}(y) \stackrel{def}{=} \max \left\{ \ \nu \{z: \ f(z) > y \}, \ \nu\{z: \ f(z) < -y\} \   \right\}, \ y > 0.
 $$

\  It is known that  by virtue of Tchebychev - Markov inequality: if  $  \ \zeta \ne 0,  $ and $ \  \zeta \in G\psi(a,b), \ $ then

  $$
  T_{\zeta} ( y) \le \exp \left( \ - h_{\psi}^* (\ln ( y/||\zeta||) ) \ \right), \ y \ge ||\zeta||, \eqno(2.2)
  $$
where

$$
h(p) = h[\psi](p)  \stackrel{def}{=} p \ \ln \psi(p), \  a \le p < b;
$$
and $ \ h^*(\cdot) \ $ denotes a famous Young-Fenchel, or Legendre transform for the  function $ \ h(\cdot): \ $

$$
h^*(v) \stackrel{def}{=} \sup_{p \in (a,b)} (pv - h(p)).
$$
 \ This assertion is alike to the famous Chernoff's estimate. It allows us to deduce the exponential tails bounds for the function $ \ \zeta = \zeta(z). \ $ \par

 \vspace{4mm}

  \ Let us introduce a very popular example.  Let $ \ \gamma = \const in (0,\infty); \ $ define the following $ \ \Psi \ $ function

$$
\psi_{\gamma}(p) \stackrel{def}{=} p^{\gamma}, \ p \in [1,\infty).
$$

 \ If $ \  f \in G\psi_{\gamma}: \ ||f||G\psi_{\gamma} = K \in (0,\infty),  \  $  i.e.

$$
\sup_{p \ge 1} \left\{ \ \frac{|f|_p}{p^{\gamma}} \  \right\} = K < \infty,
$$

 then

$$
T_f(y) \le \exp \left\{ \ - \gamma \ e^{-1} \ (y/K)^{1/\gamma} \  \right\}, \  y \ge K. \eqno(2.3)
$$

 \ For instance, the case $ \  \gamma = 1/2 \ $ correspondent with classical subgaussian functions. \par

 \ If the measure $ \ \nu \ $ is bounded, for instance  $ \ \nu(Z) = 1, \ $ then the inverse conclusion to the (2.3) is also true:
the any (measurable) function $ \ f: Z \to R \ $ satisfies the estimate (2.3), then $ \ f \in G\psi_{\gamma} \ $ and moreover
$ \ ||f||G\psi_{\gamma} = C(\gamma) \cdot K \in [0,\infty). \  $\par

 \vspace{4mm}

  \ These Grand Lebesgue Spaces (GLS) are also closely related under simple natural conditions with the so-called exponential Orlicz ones. Namely,
introduce the following {\it exponential} Young-Orlicz function

$$
N_{\psi}(u) = \exp \left(h_{\psi}^* (\ln |u|) \right),  \ |u| \ge 1; \ N_{\psi}(u) = C u^2, \ |u| < 1,
$$
and the correspondent Orlicz norm will be denoted by $ \ ||\cdot||L \left(N_{\psi} \right) =  ||\cdot||L (N). \ $ It was  done

$$
||\zeta||G\psi \le C_1 ||\zeta||L(N)  \le C_2 ||\zeta||G\psi, \ 0 < C_1 < C_2 < \infty. \eqno(2.4)
$$

 \  Note for instance that  for the $ \ \psi_{\gamma} \ $ function the correspondent Young-Orlicz one has a form

$$
N_{\psi_{\gamma}} (u) = \exp( \ C \ |u|^{1/\gamma} \   ), \ |u| \ge 1.
$$

 \vspace{4mm}

  \ Let us introduce the following example, with a following {\it degenerate} $ \ \psi \ - $ function. Define

$$
\psi_{(r)}(p) := 1, \ p = r, \ \psi_{(r)}(p) = + \infty, \ p \ne r. \eqno(2.5)
$$
Here $ \ r = \const \in [1,\infty). \ $ \par
 \ The classical Lebesgue - Riesz norm $ \ |f|_r \ $ coincides  with GLS one relative the $ \ \psi_{(r)}(p) \ $ function:

$$
|f|_r = ||f||G\psi_{(r)},
$$
if we take of course $ \ C/(+\infty) \ = 0. \ $  \par
 \ Thus, the classical theory of Lebesgue - Riesz spaces may be embedded onto GLS one. \par

\vspace{6mm}

 \section{Main result. Upper estimate in the GLS norm. Exactness.}

\vspace{4mm}

 \ Let us suppose that each function $ \ f_i(\cdot) \ $ belongs to some $ \  G\psi_i \  $ space:

$$
\exists (a_i, \ b_i), \ 1 \le a_i < b_i \le \infty \ \forall q_i \in (a_i, b_i) \ \Rightarrow
$$

$$
|f_i|L(q_i, X_i) \le \psi_i(q_i) \ ||f_i||G\psi_i, \ i = 1,2,\ldots, d. \eqno(3.0)
$$

 \ Of course, each  these $ \ \Psi \ $ function $ \psi_i(q_i) \ $  be choosed as a natural ones:

$$
\psi_i(q_i) := |f_i|_{q_i},
$$
if they are finite  still for some values $ \ q_i \in (1,\infty); \ $  then they are finite inside certain non-trivial interval
$ \ (a_i, \ b_i); \ 1 \le a_i < b_i \le \infty. \ $\par

 \ Define

$$
G(\vec{q}) \stackrel{def}{=} \overline{K}(\vec{q})  \cdot \prod_{i=1}^d \psi_i^{\alpha_i}(q_i), \eqno(3.1)
$$
and
$$
F_{\vec{\alpha}}[\vec{f}](\vec{q}) \stackrel{def}{=} | \ \vec{f} \ |_{\vec{\tau}(\vec{q})}^{ \vec{\alpha} }. \eqno(3.2)
$$

 \ We will proceed from the obtained before estimate

$$
|g|_p = |V(f_1, f_2, \ldots, f_d)|_{\Theta(\vec{q})} \le \overline{K}(\vec{q}) \cdot  | \ \vec{f} \ |_{\vec{\tau}(\vec{q})}^{ \vec{\alpha} }. \eqno(3.3)
$$

 \ We derive after substituting

$$
|g|_p \le G(\vec{q})  \cdot F_{\vec{\alpha}}[\vec{f}](\vec{q}),
$$
if of course $ \ p = \Theta(\vec{q}). \ $ \par

 \ Let us introduce the set of "layers"

$$
R(p) \stackrel{def}{=} \{ \ \vec{q}; \vec{q} \in D; \ \Theta(\vec{q}) = p \ \}, \ p \in(a,b), \ 1 \le a < b \le \infty.\eqno(3.4)
$$
 and define

$$
  \kappa(p) = \kappa[\vec{f}, V(\cdot)](p) \stackrel{def}{=} \inf_{\vec{q} \in R(p)} \left[ G(\vec{q}) \ \cdot F_{\vec{\alpha}}[\vec{f}](\vec{q}) \ \right]. \eqno(3.5)
$$

 \ We proved really the following main result of this report. \par

\vspace{6mm}

{\bf Theorem 3.1.} We assert in fact that under formulated above restrictions and notations

\vspace{4mm}

$$
||g[\vec{f}, V]||G\kappa \le 1, \eqno(3.6)
$$
with correspondent exponential tail estimation. \par

\vspace{5mm}

 \ {\it Let us discuss now the exactness of the estimate of theorem 3.1.} It is true still in the so-called "one - dimensional case"  $ \ d = 1, \ $ see [32], [33], [34]. \par

 \ Note that the exactness of our estimates holds true if for instance the each function $ \psi_i(p) \ $ coincides  correspondingly  with natural function for
the function  $  \ f_i: \ \psi_i(p) = |f_i|_p. \ $ \par

\vspace{4mm}

 \ The multivariate case $ \ d \ge 2 \ $  may be investigated quite analogously.  In detail,  denote alike in [33]

$$
U(\psi,f) = \left[ \ \frac{ ||V[f]|| G\psi}{||g||G\kappa} \ \right];
$$
and

$$
\overline{U} = \sup_{\psi \in \Psi} \sup_{0 \ne f \in G\psi} U(\psi,f); \eqno(3.7)
$$
then

$$
 \overline{U} = 1.\eqno(3.8)
$$

  \vspace{5mm}

 \section{Examples.}

 \vspace{4mm}

 \ We will use the following auxiliary facts. \par

 \vspace{4mm}

 {\bf Lemma 4.1.}

$$
\min_{\alpha,\beta}  \left[ \ \ \alpha^{\gamma_1} \ \beta^{\gamma_2}: \ \alpha,\beta > 0, \ 1/\alpha + 1/\beta = 1 \  \right] =
\frac {(\gamma_1 + \gamma_2)^{\gamma_1 + \gamma_2} }{\gamma_1^{\gamma_1} \  \gamma_2^{\gamma_2}}.
$$
 \ Here $ \ \gamma_1, \ \gamma_2 = \const > 0.   \  $\par

 \vspace{4mm}

 {\bf Lemma 4.2.}

$$
\min_{p_1, p_2 \ge 1} \left[ \ p_1^{\gamma_1}  p_2^{\gamma_2}: \ p_1^{-1} + p_2^{-1} = p^{-1} \ \right] = p^{\gamma_1 + \gamma_2} \cdot
\frac {(\gamma_1 + \gamma_2)^{\gamma_1 + \gamma_2} }{\gamma_1^{\gamma_1} \  \gamma_2^{\gamma_2}}.
$$
 \ As  before, $ \ \gamma_1, \ \gamma_2 = \const > 0.   \  $\par

 \vspace{4mm}

 \ We return to the considered examples (1 - 10). In each cases we assume that

 $$
 f_i \in G\psi_{\gamma_i}, \ \exists \gamma_i \in (0, \infty), \eqno(4.0)
 $$
and that $ \ d \ge 2; \ $ the case $ \ d = 1 \ $  is investigated, e.g. in [32], [33].  \par

\vspace{4mm}

{\bf Example 4.1.} Suppose $ \ f_1(\cdot) \in G\psi_{\gamma_1},  \ f_2(\cdot) \in G\psi_{\gamma_2}, \ $ and as in the example 1 $ \  g(x) = f_1(x) \cdot f_2(x).  \  $
Denote

$$
\theta := \gamma_1 + \gamma_2.
$$

 \ We deduce by virtue of Lemma 4.1 that $ \ g(\cdot) \in G\psi_{\theta} \ $ and herewith

$$
||g||G\psi_{\theta}   \le \frac {(\gamma_1 + \gamma_2)^{\gamma_1 + \gamma_2} }{\gamma_1^{\gamma_1} \  \gamma_2^{\gamma_2}} \cdot
 ||f_1||G\psi_{\gamma_1}  \ ||f_2||G\psi_{\gamma_2}. \eqno(4.1)
$$

\vspace{4mm}

{\bf Example 4.2. Tensor product.} Suppose that $ \ f_1(\cdot) \in G\psi^{\gamma_1},  \ f_2(\cdot) \in G\psi^{\gamma_2}, \ $ and  that as in the example 2
$ \  g(x_1, x_2) = f_1(x_1) \cdot f_2(x_2), \ x_1 \in X_1, \ x_2 \in X_2.  \  $ We apply the relations (1.5), \ (1.5a) and deduce again $ \ g(\cdot) \in G\psi_{\theta} \ $
and furthermore

$$
||g||G\psi_{\theta}  \le ||f_1||G\psi_{\gamma_1}  \ ||f_2||G\psi_{\gamma_2} \eqno(4.2)
$$
with the exact value of the  coefficient "1." \par

 \ More generally, if for certain $ \ \Psi \ -  \ $ functions $ \ \nu_1 = \nu_1(p), \ \nu_2 = \nu_2(p) \   \Rightarrow f_i \in G\nu_i, \ i = 1,2, \ $  then

 $$
 ||g||G(\nu_1 \cdot \nu_2) \le ||f_1||G\nu_1 \cdot ||f_2||G\nu_2. \eqno(4.2a)
 $$

\vspace{4mm}

 \ Note that the equality in the last estimate (4.2a) may be attained if for instance both the $ \ \Psi \ - \ $  functions are correspondingly the natural functions for $ \ f_1, \ f_2:\ $

$$
\nu_1(p) = |f_1|_p,  \ \nu_2(p) = |f_2|_p.
$$

\vspace{4mm}

{\bf Example 4.3.  Integral bilinear operator.} \par

\vspace{4mm}

 \ We return  to the  integral bilinear (regular) operator (1.6). Suppose for simplicity that $ \  \mu(X) = \mu_j(X_j) = 1 \ $ and that

$$
\overline{L} := \vraisup_{x,x_1,x_2} |L(x, x_1,x_2) | < \infty;
$$
then

$$
|g|_p \le \overline{L} \cdot |f_1|_{p} \ |f_2|_{p}.
$$

 \ One can repeat the previous considerations: if for certain $ \ \Psi \ -  \ $ functions $ \ \nu_1 = \nu_1(p), \ \nu_2 = \nu_2(p) \   \Rightarrow f_i \in G\nu_i, \ i = 1,2, \ $  then

 $$
 ||g||G(\nu_1 \cdot \nu_2) \le  \overline{L} \cdot ||f_1||G\nu_1 \cdot ||f_2||G\nu_2. \eqno(4.3)
 $$

\vspace{4mm}

{\bf Example 4.4.  Classical convolution.} \par

\vspace{4mm}

 \ Let $    \  f_1 \in G\psi_{\gamma_1}, \   f_2 \in G\psi_{\gamma_2}, \  X = X_1 = X_2 = R^n   \ $  or  $ \  X = X_1 = X_2 = [0,2\pi]^n \ $ and
$ \  g = f_1*f_2. \ $  One can apply the proposition of Lemma 4.2 and Beckner's estimate:

$$
||g||G\psi_{\theta} \le \frac {(\gamma_1 + \gamma_2)^{\gamma_1 + \gamma_2} }{\gamma_1^{\gamma_1} \  \gamma_2^{\gamma_2}} \ ||f_1||G\psi_{\gamma_1}  \ ||f_2||G\psi_{\gamma_2}. \eqno(4.4)
$$
 \ A slight generalization: let $ \ f_1 \in G\zeta_1, \ f_2 \in G\zeta_2 \ $  for some $ \ \Psi \ -  \ $ functions $ \ \zeta_1, \ \zeta_2; \ $ define a new such a function

$$
\zeta(p) = \zeta[\zeta_1, \ \zeta_2](p) :=
$$

$$
\inf \left\{G(p, \ p_1, \ p_2) \ \zeta_1(p_1) \ \zeta_2(p_2): \  \ p_1, \ p_2 \in (1,\infty), \ 1/p_1 + 1/p_2 =  1 + 1/p \ \right\}.
$$
 \ We deduce a non - refined up to multiplicative constant  in general case convolution  estimate in the GLS terms

$$
||g||G\zeta \le ||f_1||G\zeta_1 \ ||f_2||G\zeta_2. \eqno(4.4a)
$$

 \  Recall that in this example

$$
1 \le \frac{1}{p_1} + \frac{1}{p_2} \le 2.
$$

\vspace{5mm}

{\bf Example 4.5.  Infimal convolution, see example 5.} \par

\vspace{4mm}

 \  A very  simple estimate:

$$
|g(x)| = |f_1 \ \Box \ f_2 |(x) \le |f_1(x/2)| + |f_2(x/2)|, \ x \in R^d. \eqno(4.5)
$$

 \ Further, suppose $ \  f_1,  f_2 \in G\psi \  $ for some $ \ \psi \in \Psi. \ $ It follows immediately from (1.8a) and (1.8b) that

$$
|g|_p \le 2^{d/p} \ \psi(p) \ \left[ \ ||f_1||G\psi + ||f_2||G\psi \ \right]  \le 2^d \ \psi(p) \  \left[ \ ||f_1||G\psi + ||f_2||G\psi \ \right].
$$
 \ Thus, $ \ g(\cdot) \in G\psi  \  $ and herewith

$$
||g||G\psi \le 2^d \left[ \ ||f_1||G\psi + ||f_2||G\psi \ \right]. \eqno(4.5a)
$$

\vspace{5mm}

{\bf Examples 4.6, \ 4.7. Pseudo - differential product and Hilbert operation.} \par

\vspace{4mm}

 \ These cases may be investigated quite analogously to the one for ordinary convolution operation and may be omitted. \par

\vspace{5mm}

{\bf Example 4.8. Maximal operator, see (1.11).} \par

\vspace{4mm}

 \ Suppose

$$
\exists \gamma \in (0,\infty) \ \Rightarrow f_i \in G\psi_{\gamma}, \eqno(4.8a)
$$
and recall that here $ \ d \ge 2. \ $ \par
 \ We deduce solving  the following extremal problem

$$
\prod_{i=1}^d \left\{ \frac{p_i^{\gamma +1}}{p_i - 1}   \right\} \to \min
$$
subject to the limitation

$$
l: \  \sum_{i=1}^d \frac{1}{p_i} = \frac{1}{p}:
$$

$$
Z := \min_{(l)} \prod_{i=1}^d \left\{ \frac{p_i^{\gamma +1}}{p_i - 1}   \right\} = \left( \ \frac{d p}{d p - 1} \  \right)^{d(\gamma + 1)} \le
C_1(d)\ p^{d(\gamma + 1)},
$$
so that
$$
|g|_p \le C_2(d) \ p^{d(\gamma + 1)} \ \prod_{i=1}^d ||f_i||G\psi_{\gamma}, \eqno(4.8b)
$$
or equally

$$
||g||G\psi_{d(\gamma + 1)} \le C_2(d) \prod_{i=1}^d ||f_i||G\psi_{\gamma}. \eqno(4.8)
$$

\vspace{5mm}

{\bf Example 4.9. Hausdorff's  operation.}

\vspace{5mm}

 \ We find repeating the considerations and notations of the previous example that under conditions of the example 9 for the function of the form

 $$
g(x) = H_{\Phi, \vec{A} } [\vec{f}](x) := \int_{R^n} \frac{\Phi(t)}{|t|^n} \ \prod_{i=1}^m f_i(A_i(t)) \ dt, \ x \in R^n
 $$
we have

$$
|g|_p \le C_3(d, \Phi, \vec{A}) \ p^{d(\gamma + 2)} \ \prod_{i=1}^d ||f_i||G\psi_{\gamma},
$$
or equally

$$
||g||G\psi_{d(\gamma + 2)} \le C_3(d, \Phi, \vec{A}) \prod_{i=1}^d ||f_i||G\psi_{\gamma}. \eqno(4.9)
$$

\vspace{5mm}

\vspace{4mm}

{\bf Example 4.10. Toeplitz operators on sequence spaces.} See 1.13. \par

\vspace{4mm}

 \ Suppose as above that $ \ f \in G\psi_{\gamma_1}, \ x \in G\psi_{\gamma_2}, \ \gamma_1, \ \gamma_2 = \const > 0, \ $ and that

$$
\frac{1}{p_1} + \frac{1}{p_2} = \frac{1}{p'}, \  p' := \frac{p}{p-1} > 1.
$$

 \ Introduce a new $ \ \Psi \ - \ $ function

$$
\tau_{\gamma_1, \ \gamma_2}(p) \stackrel{def}{=} \frac{(\gamma_1 + \gamma_2)^{\gamma_1 + \gamma_2}}{\gamma_1^{\gamma_1} \ \gamma_2^{\gamma_2}} \cdot
\left\{ \frac{p}{p-1}  \right\}^{\gamma_1 + \gamma_2}.
$$

 \ We conclude on the basis of Lemma 4.2 and theorem 1

$$
||g||G\tau_{\gamma_1, \ \gamma_2} \le ||f||G\psi_1 \ \cdot ||x||G\psi_2. \eqno(4.10)
$$

  \vspace{5mm}

 \section{Tail description of obtained results.}

 \vspace{4mm}

 \ Many of obtained estimates may be expressed in the terms of tail behavior of used functions. Namely, suppose for simplicity in this section that all our measures are bounded:
  $ \ \mu_i(X_i) = \mu(X) = 1. \ $\par

\vspace{4mm}

 \ Let us return  at first  to the example 4.1, i.e. to the multiplicative operation;  assume that

$$
T_{f_i}(u) \le \exp \left[ \ -  (u/S_i)^{1/\gamma_i} \ \right], \ u \ge 1, \ S_i = \const > 0, \
$$

$$
i = 1,2, \ldots, d, \ \gamma_i \in (0,\infty); \eqno(5.0)
$$
then

$$
T_{f_1 \ f_2} (u) \le \exp \left\{ - C(\gamma_1, \ \gamma_2) \ \left[ \ u/(S_1 \ S_2)\right]^{ \gamma_1 \ \gamma_2/(\gamma_1 + \gamma_2)  } \  \right\}, u \ge 1. \eqno(5.1)
$$

 \vspace{4mm}

 \ At the same estimation (5.1) holds true also in the examples (4.2), (4.3) and (4.4), of course with another constants instead $ \  C(\gamma_1, \ \gamma_2). \ $ \par

 \vspace{4mm}

 \ Let us pay our attention to the example 4.8, devoting to the maximal operations. The assumption (4.8a) may be rewritten in the case when $ \ \mu_i(X_i) = 1 \ $ as follows

$$
\exists \gamma \in (0,\infty), \ \exists k_i \in (0,\infty) \ \Rightarrow T_{f_i}(u) \le \exp( \ - (u/k_i)^{ 1/\gamma  } \ ), \ i = 1,2,\ldots, d; \eqno(5.2)
$$
and proposition  of this  example - as follows

$$
T_g(u) \le \exp \left\{ \ - C_2(\vec{k}, \gamma_1, \gamma_2) \ u^{1/d(\gamma + 1)} \ \right\}, \ u \ge 1. \eqno(5.3)
$$

\vspace{4mm}

 \ Alike result holds true for the Haussdorf operation, see example 4.9. Indeed, impose on the functions $ \ \{  \ f_i \ \} \ $ again the condition  (5.2)
  as well as the conditions of boundedness of our measures $ \ \mu_i, \ \mu. \ $
Then

$$
T_g(u) \le \exp \left\{ \ - C_2(\vec{k}, \gamma_1, \gamma_2) \ u^{1/d(\gamma + 2)} \ \right\}. \eqno(5.3)
$$

 \vspace{5mm}

 \section{Concluding remarks.}

 \vspace{4mm}

 \ \hspace{5mm} {\bf A.} It is interest by our opinion  to investigate the feature of  compactness  of these multivariate operations, as well as to obtain the exact value of appeared constants
in the Grand Lebesgue Spaces setting. \par

\vspace{4mm}

 \ {\bf B.} Analogous research may be provided for another {\it singular} multivariate operations,  say, for the bilinear fractional Riesz operation [19]

$$
R[f_1,f_2](x) = \int_{R^{2n}} \frac{f_1(y) \ f_2(z) \ dy \ dz}{ \ \left( \ |x-y|^2 + |x-z|^2 \  \right)^{n-\beta}  \ }, \ \beta \in (0,n)
$$
as well as oscillator multilinear integral operations [17];  commutators of singular integral operations [17], [40]; bilinear Bochner-Riesz means [24], [28], and so one. \par

\vspace{4mm}

\ {\bf C.} One of interest application of the estimates for  multivariate  operation to a system of quadratic derivative nonlinear Schr\"odinger equations is represented in a recent article [21]. \par

 \vspace{5mm}

 \section{Acknowledgement.}

 \vspace{4mm}

 \ Authors are grateful to prof. E. Liflyand (Bar-Ilan University, Israel) for the statement of the considered in this report problem and fruitful discussions. \par

\vspace{5mm}

\begin{center}

 {\bf References.}

\end{center}

\vspace{4mm}

{\bf 1. J.A.Barrionevo, Loucas Grafacos, Danqing He, Petr Honzik, and Lucas Oliveira. } {\it  Bilinear spherical maximal function. } \\
arXiv:1704.03586v1 [math.CA] 12 Apr 2017 \\

\vspace{4mm}

{\bf 2. W.Beckner.} {\it Inequalities in Fourier analysis.} Annals of Math., 102, (1975), 159-182.\\

\vspace{4mm}

{\bf 3. Benedek A. and Panzone R.} {\it The space $ \ L(p) \ $ with mixed norm.}
Duke Math. J., 28, (1961), 301-324.\\

\vspace{4mm}

{\bf 4. Bernicot  Fred´eric.} {\it Local estimates and global continuities in Lebesgue spaces for bilinear operators.}\\
arXiv:0801.4088v3 [math.FA] 8 Feb 2008 \\

\vspace{4mm}

{\bf 5. Besov  O.V.,  Ilin  V.P.,  Nikolskii  S.M.} {\it Integral representation of functions and imbedding theorems.}
Vol.1; Scripta Series in Math., V.H.Winston and Sons, (1979), New York, Toronto, Ontario, London.\\

\vspace{4mm}

{\bf 6. H.J.Brascamp and E.H.Lieb.} {\it Best constants in Young's inequality, its converse and its generalization to more than three functions.}
Adv. Math. 20 (1976), 151 \ - \ 173. \\

\vspace{4mm}

{\bf 7. G.Brown, F.M´oricz.} {\it Multivariate Hausdorff operators on the spaces } $ \ L_p(R^n), \ $ J. Math. Anal. Appl. 271 (2002), 443-454. \\

\vspace{4mm}

{\bf 8. Buldygin V.V., Kozachenko Yu.V.} {\it About subgaussian random variables.}
Ukrainian Math. Journal, 1980, 32, No 6, 723-730, (in Russian).

\vspace{4mm}

{\bf 9. Buldygin V.V., Kozachenko Yu.V.} {\it Metric Characterization of Random Variables and Random Processes.} 1998, Translations of Mathematics
Monograph, AMS, v.188..

\vspace{4mm}

{\bf 10. Capone C., Fiorenza A., Krbec M.} {\it  On the Extrapolation Blowups in the Lp Scale.} Manuscripta Math., 99(4), 1999, p. 485-507. \\

\vspace{4mm}

{\bf 11. Wei Chena,  Chunxiang Zhub.} {\it Weighted estimates for the multilinear maximal function on the upper half-spaces.} \\
arXiv:1705.04939v1 [math.AP] 14 May 2017 \\

\vspace{4mm}

{\bf 12. Nguen Minh Chuong, Dao Van Duong.} {\it Multilinear Haussdorf operators on some function spaces with variable exponent. } \\
arXiv:1709.08185v1 [math.CA] 24 Sep 2017 \\

\vspace{4mm}

{\bf 13. Wet Dai and Guozhen Lu.} {\it L(p) Estimates for the bilinear Hilbert transform for $  \ 1/1 < p < 2/3: \ $ a counterexample and generalization to
non - smooth symbol.} \\ arXiv:1409.3875v2 [math.CA] 24 Oct 2014 \\

\vspace{4mm}

{\bf 14. Wang Ding-huai and Zhou Jiang.}  {\it  Necessary and sufficient conditions for boundedness of commutators of bilinear Hardy-Littlewood maximal function.}
arXiv:1708.09549v1 [math.FA] 31 Aug 2017.\\

\vspace{4mm}

{\bf 15.  A.Fiorenza.} {\it Duality and reflexivity in grand Lebesgue spaces.} Collectanea Mathematica, (electronic version), 51, 2, (2000), 131-148.\\

\vspace{4mm}

{\bf 16. A. Fiorenza and G.E.Karadzhov.} {\it Grand and small Lebesgue spaces and their analogs. } 2005,
Consiglio Nationale Delle Ricerche, Instituto per le Applicazioni del Calcoto Mauro Picone, Sezione di Napoli, Rapporto tecnico n., 272/03. \\

\vspace{4mm}

{\bf 17. Maxim Gilula, Philip T.Gressman, and Lecao Xitao.} {\it Higher decay Inequalities for multilinear oscillatory integrals. }
arXiv:1612.00050v1 [math.CA] 30 Nov 2016 \\

\vspace{4mm}

{\bf 18. Loukas Grafakos, Danqing He, Lenka Slavnikov.}  \ $ \ L_2 \times L_2 \to L_1 \ $ {\it Boundedness criteria.}
arXiv:1802.09400v1  [math.CA]  26 Februar 2018\\

\vspace{4mm}

{\bf 19. Jarod Hart, Rodolfo H.Torres, and Xinfeng Wu.} {\it  Smoothing properties of bilinear operators and Leibnitz - type rules in Lebesgue and mixed Lebesgue spaces.  } \\
arXiv:1701.02631v1 [math.CA] 10 Jan 2017 \\

\vspace{4mm}

{\bf 20. Jawerth B., Milman M.} {\it Extrapolation Theory with Applications.} Mem. Amer. Math. Soc., 440, (1991).\\

\vspace{4mm}

{\bf 21. Hiroyuki Hirayama and Shinya Kinoshita.} {\it Sharp bilinear estimates and its applications to a system of quadratic derivative nonlinear Schr\"odinger equations. }
arXiv:1802.06563v1 [math.AP] 19 Feb 2018\\

\vspace{4mm}

{\bf 22. T.Iwaniec and C. Sbordone.} {\it On the integrability of the Jacobian under
minimal hypotheses.}  Arch. Rat.Mech. Anal., 119, (1992), 129-143. \\

\vspace{4mm}

{\bf 23. T.Iwaniec, P. Koskela and J. Onninen.} {\it Mapping of finite distortion:
Monotonicity and Continuity.} Invent. Math., 144, (2001), 507-531. \\

\vspace{4mm}

{\bf 24. Eunhee Jeong, Sanghyuk Lee, and Ana Vargas.} {\it  Improved bound for the bilinear Bochner-Riesz Operator. }\\
arXiv:1711.02425v1 [math.CA] 7 Nov 2017 \\

\vspace{4mm}

{\bf 25. Kozachenko Yu.V., Ostrovsky E., Sirota L.} {\it Relations between exponential tails, moments and moment generating functions
for random variables and vectors.}\\
arXiv:1701.01901v1 [math.FA] 8 Jan 2017 \\

\vspace{4mm}

{\bf 26. Kozachenko Yu. V., Ostrovsky E.I.} (1985). {\it The Banach Spaces of random Variables of subgaussian Type.} Theory of Probab. and Math. Stat. (in
Russian). Kiev, KSU, 32, 43-57. \\

\vspace{4mm}

{\bf 27. Ishwart Kunwar and Yumeng Ou.} {\it  Two - weight inequalities for multilinear commutators. }\\
arXiv:1710.07392v1 [math.CA] 20 Oct 2017 \\

\vspace{4mm}

{\bf 28. Heping Liu and Min Wang.} {\it Boundedness oh the bilinear Bochner - Riesz means in the non-Banach triangle case. } \\
arXiv:1712.09235v1 [math.FA] 26 Dec 2017\\

\vspace{4mm}

{\bf 29. Feng Liu, Qingying Xue, and K.Oz.O Yabuta.} {\it Regularity and continuity of the multilinear strong maximal operators.}\\
arXiv:1801.09828v1 [math.CA] 30 Jan 2018 \\

\vspace{4mm}

{\bf 30. Akihito Miyachi and Naohito Tomita.} {\it  Bilinear pseudo - differential operators with exotic symbols.} \\
arXiv:1801.06744v1 [math.CA] 21 Jan 2018\\

\vspace{4mm}

{\bf 31. Ostrovsky E., Sirota L.} {\it Central Limit Theorem and exponential tail estimations in mixed Grand Lebesgue Spaces.} \\
arXiv:0915.2538v1 [math.Pr] 3 Feb 2015 \\

\vspace{4mm}

{\bf 32. Ostrovsky E., Sirota L., and Rogover E.} {\it Integral Operators in bilateral Grand Lebesgue Spaces.}\\
arXiv:0912.2538v1 [math.FA] 13 Dec 2009 \\

\vspace{4mm}

{\bf 33. Ostrovsky E.} {\it Boundedness of Operators in Bilateral Grand Lebesgue Spaces with Exact and Weakly Exact Constant Calculation.}
arXiv:1104.2963v1 [math.FA] 15 Apr 2011 \\

\vspace{4mm}

{\bf 34. Ostrovsky E.I.} (1999). {\it Exponential estimations for Random Fields and its
applications,} (in Russian). Moscow-Obninsk, OINPE. \\

\vspace{4mm}

{\bf 35. Ostrovsky E.I.}  (1982). {\it Generalization of Buldygin-Kozachenko norms and
CLT in Banach space. Theory Probab. Appl., 27, V.3, p. 617-619, (in Russian). } \\

\vspace{4mm}

{\bf 36. Ostrovsky E. and Sirota L.}  {\it Multidimensional probabilistic rearrangement
invariant spaces: a new approach.} \\
arXiv:1202.3130v1 [math.PR] 14 Feb 2012 \\

\vspace{4mm}

{\bf 37. Ostrovsky E. and Sirota L.} {\it Vector rearrangement invariant Banach spaces
of random variables with exponential decreasing tails of distributions.} \\
arXiv:1510.04182v1 [math.PR] 22 Sep 2015 \\

\vspace{4mm}

{\bf 38. Thomas Str\"omberg. } { A Study of the Operation of Infimal Convolution.}
DOCTORAL THESIS 1994; 139,  DEPARTMENT OF MATHEMATICS, ISSN 0348 - 8373,
Lulea University of Technology,  ISRN HLU - TH - T - -139,  Lulea, Sweden. \\

\vspace{4mm}

{\bf 39. Nicola Thorn.} {\it Bounded multiplicative Toeplitz operators on sequence spaces.} \\
arXiv:1801.09478v1  [math.FA]  29 Jan 2018\\

\vspace{4mm}

{\bf 40. Dinghuai Wang, Jiang Zhouc and Zhidong Teng. } {\it Sharp estimates for commutators of bilinear operators on Morrey type spaces.}\\
arXiv:1703.06395v1 [math.FA] 19 Mar 2017 \\

\end{document}